\documentclass{article}

\usepackage[preprint]{neurips_2019}

\usepackage[utf8]{inputenc} 
\usepackage[T1]{fontenc}    
\usepackage{hyperref}       
\usepackage{url}            
\usepackage{booktabs}       
\usepackage{amsfonts}       
\usepackage{nicefrac}       
\usepackage{microtype}      

\usepackage{amssymb}
\usepackage{mathtools}
\usepackage{array}
\usepackage{xcolor}
\usepackage{graphicx}
\usepackage{epstopdf}
\usepackage{algpseudocode,algorithm,algorithmicx}  
\algrenewcommand\algorithmicrequire{\textbf{Precondition:}}  
\algrenewcommand\algorithmicensure{\textbf{Postcondition:}}
\usepackage{subcaption}

\usepackage{amsthm}
\newtheorem*{definition}{Definition}

\newcommand{\dd}[2]{\frac{d #1}{d #2}}
\newcommand{\pp}[2]{\frac{\partial #1}{\partial #2}}
\newcommand{\R}{\mathbb{R}}

\newcommand{\ad}{\overline{D}}

\newcommand{\av}{\overline{v}}
\newcommand{\aw}{\overline{w}}

\title{Linear Range in Gradient Descent}

\author{%
  Angxiu Ni\thanks{\url{https://math.berkeley.edu/\~niangxiu/}} \\
  Department of Mathematics, University of California, Berkeley, CA 94020, USA \\
  \texttt{niangxiu@math.berkeley.edu} \\
\and
  Chaitanya Talnikar\\
  Nvidia Corporation, Santa Clara, CA 95051, USA\\
  \texttt{chaitukca@gmail.com}\\
}

\begin{document}

\maketitle

\begin{abstract}
This paper defines linear range as the range of parameter perturbations which lead to approximately linear perturbations in the states of a network.
We compute linear range from the difference between actual perturbations in states and the tangent solution.
Linear range is a new criterion for estimating the effectivenss of gradients and thus having many possible applications.
In particular, we propose that the optimal learning rate at the initial stages of training is such that parameter changes on all minibatches are within linear range.
We demonstrate our algorithm on two shallow neural networks and a ResNet.
\end{abstract}

\section{Introduction}

Machine learning is a popular method for
approximating complex functions arising in a variety of fields like
computer vision \citep{krizhevsky2012imagenet}, speech recognition \citep{graves2013speech} and many more.
Stochastic Gradient Descent (SGD) is a common choice for optimizing
the parameters of the neural network representing the functions \citep{deeplearning_book_Goodfellow,lecun2012efficient}.
But, the correct application of the SGD algorithm requires setting an initial learning rate, or stepsize,
and a schedule that reduces the stepsize as the algorithm proceeds closer to the optimal parameter values. 
Choosing a stepsize too small results in no improvement in the cost function, 
while choosing a stepsize too large causes non-convergence.

There have been a number of approaches suggested to solve this problem. 
Adaptive SGD algorithms like AdaGrad and Adam \citep{duchi2011adaptive}
require the setting of a global stepsize and other hyperparameters.
In non-stochastic gradient descent, there are two methods popular for determining the optimal step size for a gradient. 
The first is the is trust region method \citep{trust_region_conn,trust_region_Byrd,trust_region_sorensen},
and the second is the line search method, like backtracking line search \citep{armijo1966minimization}, 
Wolfe conditions \citep{wolfe1969convergence}, and probabilistic line search \citep{mahsereci2015probabilistic}.
However, when applied to SGD, both methods solve an optimization within a prescribed region or along a direction,
which could lead to over-optimization for the current minibatch, but deterioration of other minibatches.
Moreover, trust regions methods typically use second-order models, which can be expensive to build;
line search methods typically give only upper but not lower bound on the objective change, which might lead to even more over-optimization.

The center issue of stepsize selection is the lack of a good criterion for deciding the quality of descent directions and stepsizes.
This paper provides the criterion of linear range, defined as the range of parameter perturbations having small nonlinear measurement. 
The nonlinear measurement is the relative difference between the actual state perturbations 
and the linearized state perturbations given by the tangent solution.
As an application, we propose to select stepsizes, at the initial stages of the training process, 
by imposing a `speed limit' such that all minibatches are within linear range.

The paper is organized as follows.
First, we define tangent and adjoint solutions and show their utilities and relations in sensitivity analysis.
Then, we define linear range and develop linGrad.
Finally, we demonstrate linGrad on a few networks with different architectures.

\section{Preparations}

In this section, we first define neural networks as dynamical systems.
Then in section~\ref{s:tangent},
we show that small perturbations on parameters will lead to roughly linear perturbation in all states of the network, 
which further leads to a tangent equation for linear perturbation in the objective.
Finally, in section~\ref{s:adjoint}, we show this tangent formula is equivalent to the adjoint formula,
which is a generalization of the backpropagation.
Above discussion leads to the natural conclusion that,
gradients are meaningful only when stepsizes lead to roughly linear perturbation in the states.

To start, we define a neural network with $I$ layers as a discrete dynamical system, governed by:
\begin{equation} \label{e:primal_system_diffeo}
  u_0=x
  \,,\quad
  u_{i+1} = f_i(u_i,s_i) \text{\; for \;} 0\le i \le I-1\,,
\end{equation}
where $x$ is the input, column vectors $u_i\in \R^{m_i\times 1}$ are states of neurons, 
and $s_i\in\R^{n_i \times 1}$ are parameters at the $i$-th layer to be trained,
to minimize some objective $J$, defined as:
\begin{equation} \label{e:J}
  J\left(\{u_i,s_i\}_{i=0}^I \right)= \sum_{i=0}^{I} J_i(u_i,s_i).
\end{equation}
Typically, the objective is the difference between the actual output of the network for a given input and the output specified by the data:
in this case, the objective depends only on $J_I(u_I)$.
However, for future development, we allow the objective depend on all layers of network.

\subsection{Tangent solutions} \label{s:tangent}

Assume we want to make perturbations to parameters $\{s_i\}_{i=0}^I$ in the direction of $\{\sigma_i\in\R^{n_i \times 1}\}_{i=0}^I$,
the perturbations will be  $\{ \Delta s_i=\sigma_i \psi\}_{i=0}^I$,
where $\psi\in\R$ is the stepsize, or learning rate.
When the stepsize is infinitesimal $\delta \psi$, the first order approximation at each layer is:
  $
  \delta u_0 = 0,
  \delta u_{i+1} = f_{ui} \delta u_{i} + f_{si} \sigma_i \psi,
  $
where $\delta$ indicates infinitesimal perturbations, 
and $\delta u_i$ includes perturbations propagated from previous layers.
Here $f_{ui}:= \partial f_i/ \partial u (u_i,s_i) \in \R^{m_{i+1}\times m_i}$,
and $f_{si}:= \partial f_i/ \partial s (u_i,s_i) \in \R^{m_{i+1}\times n_i}$.
There is no perturbation on $u_0$, since the input data are accurate.

Define $v_i:=\delta u_i/\delta\psi$, it is governed by the conventional inhomogeneous tangent equation:
\begin{equation} \label{e:convential inhomo tangent}
  v_0 = 0
  \,,\quad
  v_{i+1} = f_{ui} v_{i} + f_{si} \sigma_i
  \,.
\end{equation}
Later we are interested in computing $v_i\psi$,
which is easier to compute via:
\begin{equation} \label{e:vipsi}
  v_0\psi = 0
  \,,\quad
  v_{i+1}\psi = f_{ui} (v_{i}\psi) + f_{si} \Delta s_i
  \,.
\end{equation}
To extend our work to networks with architecures not satisfying equation~\eqref{e:primal_system_diffeo}, such as ResNet and neural ODE,
we need their corresponding tangent equations, which are given in appendix \ref{app:tangent for other archtec}.

Now, we can write out a tangent formula for the sensitivity $dJ/d\psi$.
More specifically, we first differentiate each term in equation~\eqref{e:J}, then apply the definition of $v_i$, we get:
\begin{equation}\label{e:djdpsi_tan}
  \dd J \psi =\sum_{i=0}^{I} \left( J_{ui} v_i + J_{si} \sigma_i\right)
\end{equation}
Here both $J_{ui}:=\partial J_i/\partial u (u_i,s_i)\in \R^{1\times m_i}$ and
$J_{si}:=\partial J_i/\partial s (u_i,s_i)\in \R^{1\times n_i}$ are row vectors.

Investigating inhomogeneous tangent solutions calls for first defining the homogeneous tangent equation: 
$ w_{i+1} = f_{ui} w_i$,
which describes the propagation of perturbation on states while the parameters are fixed.
The propagation operator $D_l^i\in \R^{m_i\times m_l}$ is defined as the matrix 
that maps a homogeneous tangent solution at $l$-th layer to $i$-th layer.
More specifically,
\begin{equation} \begin{split} \label{e:D}
  D_l^i := 
  \begin{cases}
    I_d \textnormal{ (the identity matrix)}  \,,\quad \textnormal{ when } i=l \,; \\
    f_{u,i-1} f_{u,i-2} \cdots f_{u,l+1}f_{u,l} \,, \quad \textnormal{ when } i>l  \,.
  \end{cases}
\end{split} \end{equation}

We can use Duhamel's principle to analytically write out a solution to equation~\eqref{e:convential inhomo tangent}.
Intuitively, an inhomogeneous solution can be viewed as linearly adding up homogeneous solutions,
each starting afresh at a previous layer, with initial condition given by the inhomogeneous term.
More specifically, 
\begin{equation} \label{e:v by Duhamel}
  v_0 = 0
  \,, \quad 
  v_i = \sum_{l=0}^{i-1} D^i_{l+1}f_{sl}\sigma_l \text{\; for \;} 1 \le i \le I. 
\end{equation}

\subsection{Adjoint solutions} \label{s:adjoint}

In this subsection, we first use a technique similar to backpropagation to derive an adjoint sensitivity formula,
which we then show is equivalent to the tangent sensitivity formula in equation~\eqref{e:djdpsi_tan}.

Assume we perturb the $l$-th layer by  $\delta u_l$, and let it propagate through the entire network,
then the change in the objective is:
$J+\delta J = \sum_{j=l}^{I} J_j(f_{I-1}(\cdots f_l(u_l+\delta u_l,s_l)\cdots),s_j)$.
Neglecting higher order terms, we can verify the inductive relation
$\delta J /\delta u_l =\delta J /\delta u_{l+1}  f_{ul} + J_{ul}$.
Define $\av_l: = \delta J /\delta u_{l} \in \R^{1\times m_l}$,
it satisfies the conventional inhomogeneous adjoint equation:
\begin{equation} \label{e:inhomo_adjoint_diffeo}
  \av_{I+1} = 0 
  \,,  \quad 
  \av_{l} = \av_{l+1}  f_{ul} + J_{ul}
  \,.
\end{equation}
Notice the reversed order of layers.
Here, the terminal condition is used because we can assume there is $(I+1)$-th layer which $J$ does not depend on.
Hence, the adjoint sensitivity formula is:
\begin{equation} \begin{split} \label{e:adjoint sensitivity directly}
  \dd J \psi = \sum_{l=0}^I \frac{\delta J}{\delta u_l} \pp {u_l}\psi + J_{sl}\sigma_l
  = \sum_{l=1}^I \av_l f_{s,l-1}\sigma_{l-1} + \sum_{l=0}^I J_{sl}\sigma_l \,,
\end{split} \end{equation}
where $\partial u_0 /\partial \psi =0 $ as $u_0$ is fixed, and $ \partial u_l /\partial \psi = f_{s,l-1}\sigma_{l-1}$.
Notice that $\partial u_l /\partial \psi$ is not the tangent solution $ v_l= \delta u_l /\delta \psi$,
since $\delta u_l$ in the definition of tangent solution includes not only perturbation due to change in $s_{l-1}$,
but also the perturbation  propagated from the previous layer.
In other words, in the tangent formula, the propagation of perturbations on states is included in $v_l$,
whereas in the adjoint formula such propagation is included in $\av_l$.

The advantage of the adjoint sensitivity formula, comparing to the tangent formula,
is a clearer view of how the sensitivity depends on $\sigma_i$,
which further enables us to select the direction for perturbing parameters,  $\{\sigma_i\}_{i=0}^I$.
Not surprisingly, the inhomogeneous adjoint solution is a generalization of the backpropagation.
To illustrate this, if we set the objective to take the common form, $J = J_I(u_I)$,
then $\av_{l} = J_{uI} f_{u, I-1} \cdots f_{u,l}$.
The gradient of $J$ to parameters, given by the backpropagation, is:
\begin{equation} \begin{split}
  \partial J / {\partial s_l} &= J_{uI} f_{u, I-1} \cdots f_{u,l+1} f_{sl} = \av_{l+1}f_{sl} \,.
\end{split} \end{equation}

The sensitivity can be given by either a tangent formula in equation~\eqref{e:djdpsi_tan},
or an adjoint formula in equation~\eqref{e:adjoint sensitivity directly},
hence, the two formula should be equivalent.
Since later development heavily depends on this equivalence, we also prove it directly.
To start, first define the homogeneous adjoint equation:
$ \aw_{l} = \aw_{l+1} f_{ul}$, 
where $\aw_l \in \R^{1\times m_i}$ is a row vector.
The adjoint propagation operator $\ad_i^l$ is the matrix which, multiplying on the right of a row vector,
maps a homogeneous adjoint solution at $i$th layer to $l$th layer.
A direct computation shows that $D_l^i = \ad_i^l$.
Using Duhamel's principle with reversed order of layers, 
we can analytically write out the inhomogeneous adjoint solution:
\begin{equation} \label{e:av by Duhamel}
  \av_{I+1} = 0 
  \,,  \quad 
  \av_l = \sum_{i=l}^{I} J_{ui} \ad^l_{i}
  \text{\; for \;} 0\le l\le I 
  \,.
\end{equation}
To directly show the equivalence between tangent and adjoint formula, 
first substitute equation~\eqref{e:v by Duhamel} into \eqref{e:djdpsi_tan},
change the order of the double summation, then assemble terms with the same $f_{sl}\sigma_l$:
\begin{equation} \begin{split} \label{e:djdpsi_adj}
  \dd J \psi 
  &=\sum_{i=1}^{I} J_{ui} v_i + \sum_{i=0}^{I} J_{si} \sigma_i
  =\sum_{i=1}^{I} \sum_{l=0}^{i-1} J_{ui} \ad_i^{l+1} f_{sl}\sigma_l + \sum_{i=0}^{I} J_{si} \sigma_i \\
  &=\sum_{l=0}^{I-1} \left( \sum_{i=l+1}^{I} J_{ui} \ad_i^{l+1} \right) f_{sl}\sigma_l + \sum_{i=0}^{I} J_{si} \sigma_i
  =\sum_{l=0}^{I-1}  \av_{l+1} f_{sl}\sigma_l + \sum_{i=0}^{I} J_{si} \sigma_i \,.
\end{split} \end{equation}

\section{Linear range}

\subsection{Definition}

Assuming that the direction to perturb parameters, $\{\sigma_i\}_{i=0}^{I}$, has been decided,
we still need to specify the stepsize $\psi$ to get the new parameters, 
the selection of which is the topic of this section.
There are two equivalent methods for computing the sensitivity $dJ/d\psi$:
the tangent formula in equation~\eqref{e:djdpsi_tan} and the adjoint formula in equation~\eqref{e:djdpsi_adj}.
The adjoint formula is useful for deciding $\sigma_i$, 
and the tangent formula is useful for checking the effectiveness of sensitivity as the tangent solution has the ability to predict the linear change in objective after slightly perturbing the parameters.
A sufficient condition for the approximate linearity,
is that the perturbation in all of the states are roughly linear to $\{ \Delta s_i\}_{i=0}^I$.

To elaborate, we first define a nonlinear measurement for the perturbations in the states of one network,
\begin{equation}  \label{e:nonlinear measurement}
  \varepsilon := \frac 1{I} \sum_{i=1}^I 
  \frac{
  \|u_{new,i} - u_{old,i} - v_i \psi\| 
  }
  {
  \|v_i\psi\| 
  }\,,
\end{equation}
where $v$ is the conventional tangent solution, 
$u_{old}$ and $u_{new}$ are the states before and after parameter change,
subscript $i$ indicates the layer,
and the norm is $l^2$.
Assume that we can use Taylor expansion for $u_{new}$ around $u_{old}$,
and that $v = \delta u_{new}/ \delta \psi$ is non-zero,
we have $u_{new} = u_{old} + v\psi + v'\psi^2 + O(\psi^3)$, where $v'$ is some unknown constant vector.
Hence $ \varepsilon = C \psi + O(\psi^2) $ for some constant $C$,
and for small $\psi$ we may regard $\varepsilon$ as linear to $\psi$.
With the above description of nonlinear measurement,
we can finally define linear range.

\begin{definition}
  Given a network and an input data,
  the $\varepsilon^*$-linear range on parameters is the range of parameter perturbations such that
  $\varepsilon \le \varepsilon^*$.
  The linear range on objective and on states are the image sets of the linear range on parameters.
\end{definition}

\subsection{Gradient descent by linear range} \label{s:linGrad}

Linear range is a criterion that can be used in many ways.
In this subsection, we use it to develop linGrad, which is a stochastic gradient descent (SGD) method.
In linGrad,
the stepsize is determined by a subset of samples in all minibatches, 
such that the perturbation on parameters are just within the $\varepsilon^*$-linear range.
More specifically, for each one out of several minibatches, we use the current $\psi$ to compute $\varepsilon$.
Since $\varepsilon$ is linear to $\psi$ when $\psi$ is small, 
$\psi^*=\psi \varepsilon^* /\varepsilon$ is the $\varepsilon^*$-linear range on stepsize for this minibatch.
We update $\psi$ to be the smallest $\psi^*$ within a finite history.

Algorithm~\ref{a:lingrad} lists steps of linGrad.
We suggest to use $0.3\le \varepsilon^* \le 1$,
so that the stepsize is not too small, yet the gradient is still meaningful.
Our experiments show that above range of $\varepsilon^*$ yields smaller than 10 times difference in stepsizes,
meaning that the linear range criterion reduces the possible range of optimal stepsizes to within an order of magnitude.

$N_{hist}$ should be chosen by statistical significance, for example $N_{hist}\ge 50$,
such that the max of $\psi^*$ over sampling minibatches is approximately the true max over all minibatches.
We also require $N_{hist}N_{lin}\le C N_b$ for some $C$ of order $O(1)$, 
so that only recent linear ranges affects the selection of current stepsize.
In fact, we found in our experiments that linGrad is robust to the selection of $N_{lin}$ and $N_{hist}$ once above conditions are satisfied.

\begin{algorithm}
  \caption{linGrad: Linear range gradient descent (with fixed $\varepsilon^*$)}  \label{a:lingrad}
\begin{algorithmic}[1]
  \Require $\varepsilon^*$; empty list $L$; $N_s$ samples in a minibatch;  $N_b$ minibatches; $N_{hist}$; $N_{lin}$.
  \For {each epoch,}
  \For {each minibatch,}
    \For{$n\gets 1, N_s$} 
    \Comment{The subscript $n$ is omitted sometimes.}
      \State Compute $u_{old}$ using parameters $\{s_i\}_{i=0}^I$ and input data.
      \State Compute adjoint solution $\{\av_l\}_{l=1}^I$.
      \State Select $\{\sigma_i\}_{i=0}^{I}$ according to predetermined rules.
      \State Compute tangent solution $\{v_i\psi\}_{i=1}^I$ by equation~\eqref{e:vipsi}. \label{l:extra start}
      \State Compute new states $u_{new}$ using parameters $\{s_i + \sigma_i\psi\}_{i=0}^I$.
      \State Compute $\varepsilon_n$ for this sample using equation~\eqref{e:nonlinear measurement}.
    \EndFor
    \State Compute $\varepsilon = (\sum_{n=1}^{N_s} \varepsilon_n) / N_s$. 
    \State Append $\psi^*=\psi \varepsilon^*/\varepsilon $ to the list $L$.
    \State $\psi \gets \min \{\textnormal{last } N_{hist} \textnormal{ elements in }L \}$ \label{l:extra end}
    \State Update parameters $s_i \gets s_i+\sigma \psi$. 
  \EndFor
  \Comment{Only need to perform steps \ref{l:extra start} to \ref{l:extra end} once every $N_{lin}$ minibatches.}
  \EndFor
\end{algorithmic}
\end{algorithm}

\subsection{Remarks}

Notice that the nonlinear measurement is defined over the entire network rather than just over the objective.
Since the objective is only one number, 
it may not provide adequate information for deciding where the parameter perturbations are within the linear range.
In fact, we tried defining linear measurement by objectives, and found the algorithm not robust,
for example, optimal $\varepsilon^*$ changes to settings like minibatch sizes, and for larger $\varepsilon^*$ the algorithm diverges.
We also tried adding the objective as an additional layer after the output, but still find the algorithm not robust;
further limiting the maximum contribution from the objective layer in the nonlinear measurement helps improving robustness.
We suggest readers to experiment whether and how to include objective in the definition of nonlinear measurement.

The concept of linear range is useful for other scenarios beyond linGrad.
One possible application is that it offers a criterion for comparing different descent directions:
larger linear range yields larger parameters and objective perturbations, thus faster convergence.
For example, we can use linear range to compare the gradients computed by normal backpropagation 
and by clipping gradients \citep{clip_gradients} for deep neural networks.

Another use of linGrad is to determine the initial stepsize and then change to an adaptive algorithm for stepsizes like Adam or AdaGrad.
It is also possible to increase batch size instead of decrease stepsize \citep{Byrd_big_batch,Friedlander_big_batch}.
There are also many choice in terminating criteria,
There are many choices for termination criteria for the optimization process, 
for example the optimization can be terminated when the signal-to-noise ratio, 
which is the ratio between the average and RMS of gradients, is too low \citep{De_big_batch};
or when the ratio of counter-directions, which is the pair of gradients with negative inner-products, is roughly half.
LinGrad can be added to many existing optimization algorithms and training schemes, and we suggest readers to experiment.

Although not implemented in this paper,
it is possible to obtain tangent solutions by tracing computation graphs.
The tangent propagation is a local process, just like backpropagation:
every gate in a circuit diagram can compute how perturbations in its inputs are linearly transported to the output.
Notice that here the inputs to a gate can be either states or parameters.
For cases such as convolution networks where each neuron depends only on a few neurons in the previous layer,
tangent solvers implemented using graph tracing are faster.

An easier but less accurate way to obtain tangent solutions is via finite differences.
By the definition of tangent solutions, we can see $ v_i \approx \Delta u_i / \Delta\psi$,
meaning that we can first set $\psi$ to be a small number $\delta$, say 1e-6,
then compute new states $u_{new,i}^\delta$, and then $v_i \approx (u_{new,i}^\delta-u_{old,i}) / \delta$.
This way of computing tangent solutions does not require coding a true linearized solver,
rather, it only requires running the feedforward process one more time.

\section{Applications} \label{s:applications}

\subsection{Application on an artificial data set} \label{s:DIST}

We first apply linGrad on a network where each layer is given by 
$f_i(u_i, W_i) = g(W_i u_i+ b_i)$, where $g$ is the vectorized logistic function.
Our parameters to be learned are $W_i\in\R^{m_{i+1}\times m_{i}}$ and $b_i\in\R^{m_{i+1}}$.
The perturbations on parameters are $\Sigma_i \psi = \Delta W_i$ and $\beta_i \psi = \Delta b_i$.
Our objective is defined only on the last layer as the square difference $J:=J_I(u_I) = \frac 12 \sum_{j = 1}^{m_I} (u_I^j - y^j)^2$,
where $y$ is the output data for this sample.
To adapt with our previous notations, 
we regard $(W_i, b_i)$ and $(\sigma_i, \beta_i)$ 
as one-dimensional vectors of length $n_i=m_{i+1}\times m_{i}+m_{i+1}$, obtained by flattening the matrix and appending to the vector.
Then, for programming convenience, we reshape this vector back into a matrix and a vector in the list of results below.
\begin{equation} \begin{split}
  f_{ui} = \Lambda_i W_i \,,\quad
  f_{si} \sigma_i = \Lambda_i (\Sigma_i u_i + \beta_i) \,,\quad
  J_{uI} = (u_i - y)^T \,,
\end{split} \end{equation}
where $\Lambda_i = diag[g_i(1-g_i)] \in \R^{m_{i+1}\times m_{i+1}}$
is a diagonal matrix due to differentiating the component-wise logistic function.
By either carefully managing subscripts of partial derivatives in equation~\eqref{e:vipsi} and \eqref{e:inhomo_adjoint_diffeo}, or deriving directly from the definition,
we get tangent and adjoint equations:
\begin{equation} \begin{split}
  v_0\psi = 0 
  \,,\quad  
  v_{i+1} \psi &= \Lambda_i(W_i v_i\psi + \Delta W_i u_i + \Delta b_i)  
  \,;   \\
  \av_I = J_{uI}
  \,,\quad 
  \av_{i} &= \av_{i+1} \Lambda_i W_i 
  \,.
\end{split} \end{equation}

The feedforward and backpropagation in our implementation are from the code complementing~\citep{Nielsen2015}. 
For our particular example, we use two hidden layers.
All layers have the same number of neurons, $m_i = 50$.
We first fix the network with randomly generated parameters,
then generate 50k training samples and 10k test samples by feeding this fixed network with random inputs.
Here all random numbers are from independent standard normal distribution.

For the training, initial parameters are generated randomly,
and all samples are randomly shuffled for each epoch.
We compute the nonlinear measurement and adjust stepsize every $N_{lin}=100$ minibatches,
and take stepsize as the smallest of the last $N_{hist}=\max(50,N_b/N_{lin})$ candidate values,
where $N_s$ varies.
We choose $\Sigma_i$ as the fastest descent direction:
\begin{equation}
  \Sigma_i  = - \av_{i+1}f_{si} = - \Lambda_i \av_{i+1}^T u_i^T \,, \quad
  \beta_i = - \Lambda_i \av_{i+1}^T \,.
\end{equation}

\begin{figure}[ht]
\centering
  \includegraphics[width=0.49\textwidth] {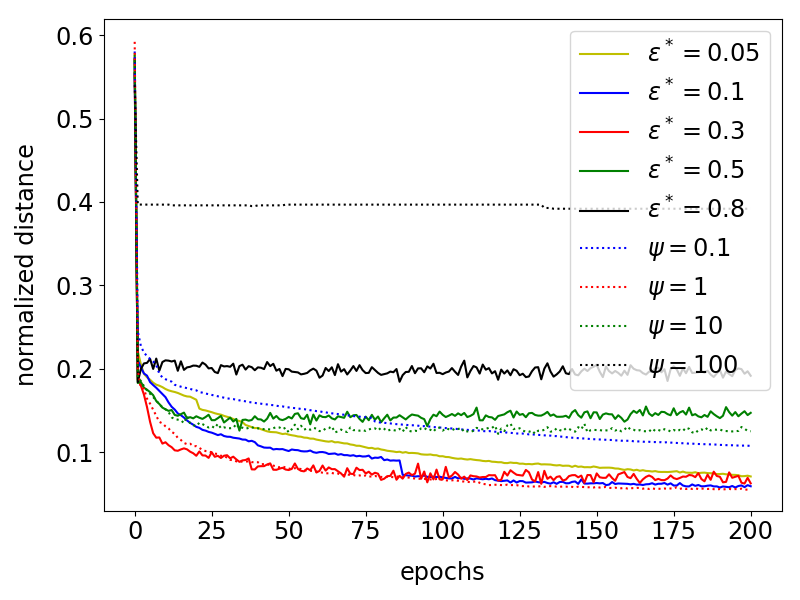}
  \includegraphics[width=0.49\textwidth] {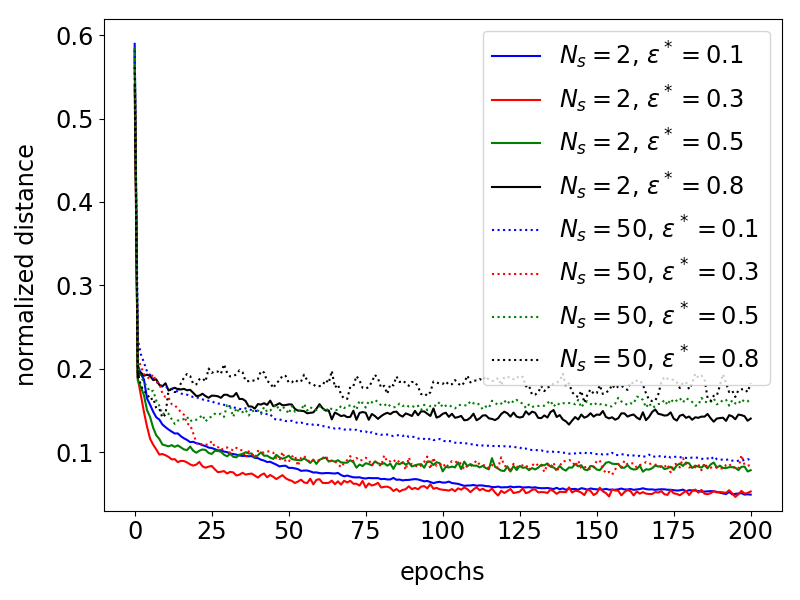}
  \caption{LinGrad applied on data generated by an artificial network.
          Each objective history is averaged over 5 runs.
          The vertical axis is average of normalized distance $(\sum_{j = 1}^{m_I} (u_I^j - y^j)^2)^{0.5}/\sqrt{m_I}$.
          Left: linGrad with minibatch size $N_s=10$ and different $\varepsilon^*$ versus SGD with different fixed stepsizes.
          Right: linGrad with different $N_s$ and $\varepsilon^*$.}
  \label{f:DIST}
\end{figure}

As we can see from the left of figure~\ref{f:DIST}, 
for batch size $N_s=10$,
comparing to SGD with fixed stepsizes,
linGrad with $\varepsilon^*=0.3$ descents the fastest, 
especially in the first 50 epochs, 
confirming that the `speed limit' during the first phase of training neural networks is given by the criterion of linear range.
In fact, if the objective function is defined as the current objective multiplied by 10,
SGD would have parameter perturbations that are 10 times larger, resulting in different convergence behavior,
whereas the linear range and hence linGrad would remain unaffected.
Moreover, from the right of figure~\ref{f:DIST}, 
we can see that $\varepsilon^*=0.3$ persists to be optimal for linGrad with different batch sizes.

\begin{figure}[ht]
\centering
  \includegraphics[width=\textwidth] {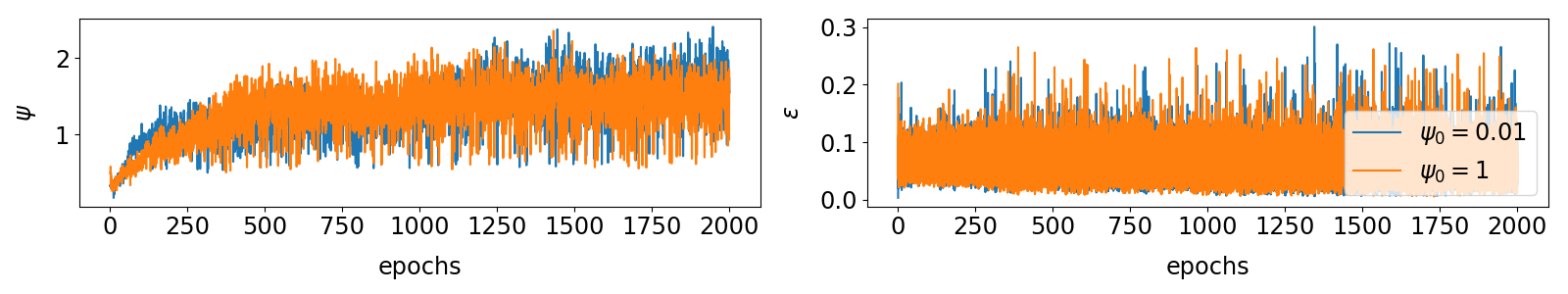}
  \caption{History of stepsize and nonlinear measurement for linGrad with $\varepsilon^*=0.3$ and $N_s=10$.}
  \label{f:DIST_psi_eps}
\end{figure}

Histories of stepsize $\psi$ and nonlinear measurement $\varepsilon$ of linGrad are shown in figure~\ref{f:DIST_psi_eps}.
We run linGrad with different initial stepsizes $\psi_0=0.01$ and $\psi_0=1$.
As shown, $\psi_0$ does not affect much:
this is expected, since $\psi_0$ is used only to infer the first $\varepsilon^*$-linear range.
This confirms that linGrad relieves the headache of choosing initial stepsizes.
Also we can see that for this shallow network, the stepsize remains at roughly the same value, 
indicating that $N_{hist}$ is statistically significant.
Finally, the nonlinear measurement remains below 0.3, 
confirming that our implementation correctly keeps the stepsize within the linear range.

\subsection{Application on MNIST}

\begin{figure}[ht]
\centering
    \includegraphics[width=0.49\textwidth] {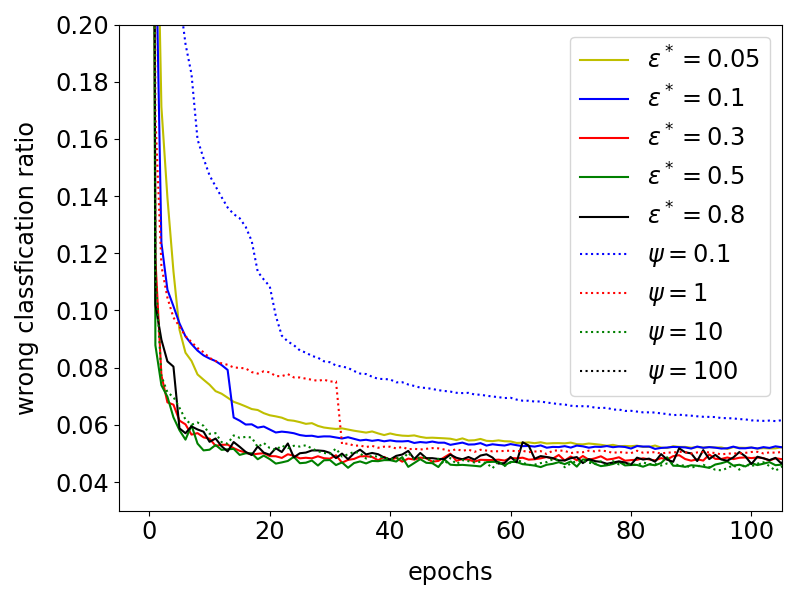}
    \includegraphics[width=0.49\textwidth] {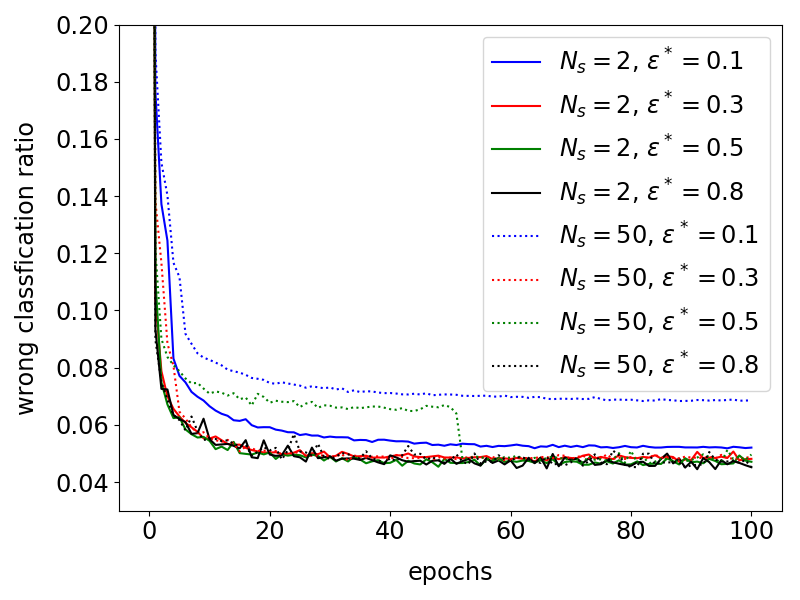}
    \caption{LinGrad applied on MNIST with the same settings as figure~\ref{f:DIST}.
      The history for SGD with $\psi=100$ does not converge and the objective value remains at 0.9, hence it is out of picture.}
  \label{f:MNIST}
\end{figure}

We then apply linGrad on MNIST with 60k training data and 10k test data.
The network has three layers, where the input layer has 764 neurons, the hidden layer 30 neurons, and the output layer 10 neurons.
The classification is done by selecting the largest component in the output layer.
Other aspects of the architecture and settings are the same as we used in section~\ref{s:DIST}.
We compare linGrad to SGD with constant stepsizes and compare linGrad with different minibatch sizes in figure~\ref{f:MNIST}.
For this problem linGrad converges fastest for either $\varepsilon^*=0.3, 0.5$ or $0.8$,
both comparable to SGD with optimal stepsize.
Again, we can see the selection of $\varepsilon^*$ is robust to $N_s$.

\subsection{Application on CIFAR-10}

Finally, we apply linGrad on the CIFAR-10 dataset using the ResNet model \citep{resnet} with 18 layers.
The size of the training dataset is 50k image samples and the number of classes for the classification task is 10. 
The size of each image is 32x32, and the total number of weights in the model is approximately 11 million.
The tangent solution of the network is computed using the finite difference method 
with a small stepsize $\delta=$1e-6 for the parameter perturbation.
The nonlinear measurement is computed using the states of the network after every residual block in ResNet.

\begin{figure}[ht]
\centering
    \includegraphics[width=0.49\textwidth] {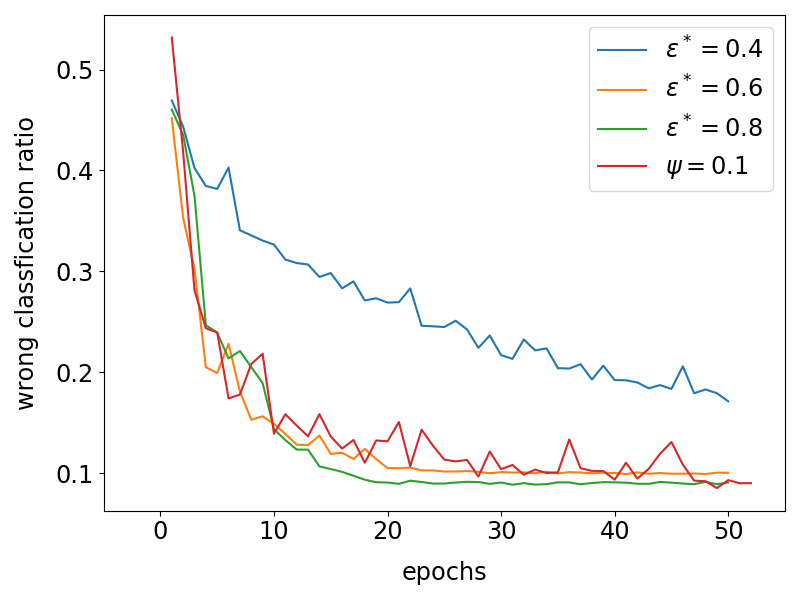}
    \includegraphics[width=0.49\textwidth] {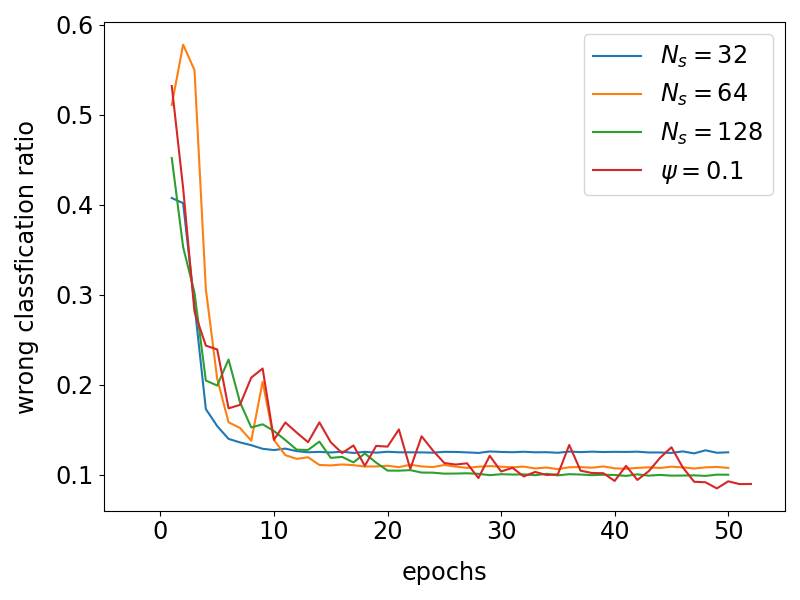}
    \caption{LinGrad applied on CIFAR-10 using ResNet-18 with $N_{lin} = 10$.
    Left: compare SGD and linGrad using $N_s = 128$.
    Right: compare linGrad with different minibatch sizes.}
  \label{f:CIFAR}
\end{figure}

We compare linGrad to SGD with constant stepsize and additionally compare linGrad with different minibatch sizes in figure~\ref{f:CIFAR} 
by computing the error on the test data.
LinGrad converges as fast as SGD for $\varepsilon^*=0.6$ and $0.8$. 
The performance of linGrad is similar across varying minibatch sizes.
Moreover, we find that for this deeper network, the stepsizes given by linGrad automatically reduces during the training.
It remains to be further investigated whether this automatic decay by linGrad fits the known optimal stepsize scheme for later stages of training.

\section{Conclusion}

This paper defines linear range and states how to compute it via comparing tangent solutions with finite differences.
Linear range is a new criteria which can be used for evaluating the quality of stepsizes and descent directions, 
and it could have many theoretical and practical applications.
In particular, we develop a stochastic gradient descent algorithm, linGrad, 
where the stepsize is given by such that all minibatches are within $\varepsilon^*$-linear range.
By applying linGrad on two shallow networks and a ResNet,
we find that the fastest convergence is obtained inside the interval $0.3\le\varepsilon^*\le 1.0$,
which corresponds to stepsize differences less than an order of magnitude.
LinGrad can be integrated with many existing gradient descent algorithms to improve the selection of stepsizes,
at least during the initial phase of the training process.

\appendix
\section{Tangent equations for other architectures} \label{app:tangent for other archtec}

To compute nonlinear measurement defined by equation~\eqref{e:nonlinear measurement} for other architectures, 
we further provide corresponding tangent equations.

For ResNet \citep{resnet}, the dynamical system corresponding to equation~\eqref{e:primal_system_diffeo}
and the inhomogeneous tangent equation corresponding to equation~\eqref{e:vipsi} are:
\begin{equation} \begin{split}
  u_0=x 
  ,\quad 
  u_{i+1} &= g\left( \sum_{j\le i} W_{j,i+1} u_j + b_i\right) 
  ;\\
  v_0\psi=0
  ,\quad 
  v_{i+1}\psi &= \Lambda_i\left( \sum_{j\le i} \left(W_{j,i+1} v_j\psi + \Delta W_{j,i+1} u_j \right)  + \Delta b_i\right)
  .\\
\end{split} \end{equation}
Further, we suggest that the objective function should not depend on the intermediate layers in a residual block,
since they have different interpretations as the input/output layers of the block.
In fact, if we do not regard intermediate layers in residual blocks as states,
then we can recover the basic form of dynamical system in equation~\eqref{e:primal_system_diffeo} with complicated $f_i$'s.

For neural ODE \citep{odenet}, the dynamical system and inhomogeneous tangent equation are:
\begin{equation} \begin{split}
  u_0 = x 
  , \quad
  \dd ut &= f(u,s)
  \,; \\
  v_0 = 0
  , \quad 
  \dd {v\psi}t &= f_u v\psi + f_s \Delta s
  \,.
\end{split} \end{equation}
Here $f_u:=\partial f/\partial u$, $f_s:=\partial f/\partial s$.
Moreover, the summation in the definition of nonlinear measurement in equation~\eqref{e:nonlinear measurement} should change to integration.

\bibliography{MyCollection}
\small {\bibliographystyle{plainnat}}

\end{document}